\documentclass[%
twocolumn,
 reprint,jcp,
 amsmath,amssymb,
 aps,superscriptaddress
floatfix,
]{revtex4-1}

\usepackage{graphicx}
\usepackage{algorithm}
\usepackage[noend]{algpseudocode}
\usepackage{longtable}
\usepackage{enumitem}
\usepackage[hidelinks]{hyperref}
\usepackage{pgfplots}
\usepackage{array}

\usepackage{colortbl,array}
\usepackage{multirow,bigdelim}
\usepackage{booktabs}

\newcolumntype{H}{>{\setbox0=\hbox\bgroup}c<{\egroup}@{}}

\newcommand{\Tau}{\mathcal{T}}
\newcommand{\mA}{\mathcal{A}}
\usepackage{alltt}
\newcommand{\bX}{{\bf X}}
\newcommand{\bR}{{\bf R}}
\newcommand{\bP}{{\bf P}}
\newcommand{\bY}{{\bf Y}}
\newcommand{\bQ}{{\bf Q}}
\newcommand{\bF}{{\bf F}}

\newcommand{\mX}{\mathcal{X}}
\newcommand{\mR}{\mathcal{R}}
\newcommand{\mP}{\mathcal{P}}
\newcommand{\mY}{\mathcal{Y}}
\newcommand{\mQ}{\mathcal{Q}}
\newcommand{\mH}{\mathcal{H}}
\newcommand{\mO}{\mathcal{O}}

\usepackage{amsthm}

\begin{document}


\title{Calculating vibrational spectra of molecules using tensor train decomposition}

\author{Maxim Rakhuba}
\email{rakhuba.m@gmail.com}
\affiliation{%
 Skolkovo Institute of Science and Technology, Skolkovo Innovation Center, Building 3, 143026 Moscow, Russia.
}
\author{Ivan Oseledets}%
 \email{ivan.oseledets@skoltech.ru}
\affiliation{%
 Skolkovo Institute of Science and Technology, Skolkovo Innovation Center, Building 3, 143026 Moscow, Russia.
}
\affiliation{Institute of Numerical Mathematics of Russian Academy of Sciences, Gubkina St. 8, 119333, Moscow, Russia.}

\date{\today}

\begin{abstract}
We propose a new algorithm for calculation of vibrational spectra of molecules using tensor train decomposition.
Under the assumption that eigenfunctions lie on a low-parametric manifold of low-rank tensors we {suggest using well-known iterative methods that utilize matrix inversion (LOBPCG, inverse iteration) and solve corresponding linear systems inexactly along this manifold.}
As an application, we accurately compute vibrational spectra (84 states) of acetonitrile molecule CH$_3$CN on a laptop in one hour using only $100$ MB of memory to represent all computed eigenfunctions.
\end{abstract}


\pacs{Valid PACS appear here}
\maketitle



\section{Introduction} \label{sec:intro}

In this work we consider time-independent Schr\"odinger equation to calculate vibrational spectra of molecules.
The goal is to find smallest eigenvalues and corresponding eigenfunctions of the Hamiltonian operator.
The key assumption that we use is that potential energy surface (PES) can be approximated by a small number of sum-of-product functions. This holds, e.g. if PES is a polynomial.

We discretize the problem using the discrete variable representation (DVR) scheme \cite{baye-dvr-1986}.
The problem is that the storage required for each eigenvector grows exponentially with dimension $d$ as $n^d$, where $n$ is number of grid points in each dimension. Even for the DVR scheme where $15$ grid points in each dimension is often sufficient to provide very accurate results we would get $\approx1$~PB of storage for a 12 dimensional problem.
This issue is often referred to as the \emph{curse of dimensionality}.
To avoid exponential growth of storage we use \emph{tensor train} (TT) decomposition \cite{osel-tt-2011} to approximate the operator and the eigenfunctions in the DVR representation, {which is known to have exponential convergence rate.
It is important to note that the TT-format is algebraically equivalent to the Matrix Product State format (MPS), which has been used for a long time in quantum information theory and solid state physics to
approximate certain wavefunctions \cite{white-dmrg-1992, ostlund-dmrg-1995} (DMRG method), see the review \cite{schollwock-2011} for more details.}


Prior research \cite{dkos-eigb-2014} has shown that the eigenfunctions often can be well approximated in the TT-format, i.e. they lie on a certain low-dimensional non-linear manifold.
The key question is how to utilize this a priori knowledge in computations.
{
We propose to use well-established iterative methods that utilize matrix inversion and solve corresponding linear systems inexactly along the manifold to accelerate the convergence.
}
Our main contributions~are:
\begin{itemize}
\item
we propose a concept of a \emph{manifold preconditioner} that explicitly utilizes information about the size of the TT-representation.
We use the manifold preconditioner for a tensor version of \emph{locally optimal block preconditioned conjugate gradient} method (LOBPCG) \cite{knyazev-lobpcg-2001}. 
We will refer to this approach as \emph{manifold-preconditioned LOBPCG} (MP LOBPCG).
The approach is first illustrated on computation of a single eigenvector (Section \ref{sec:mp1}) and then extended to the block case (Section \ref{sec:mpb}).

\item we propose tensor version of simultaneous inverse iteration (also known as {block inverse power method}), which significantly improves accuracy of the proposed MP LOBPCG.
Similarly to the manifold preconditioner the inversion is done using the a priori information that the solution belongs to a certain low-parametric manifold.
We will refer to this method as {\emph{manifold-projected simultaneous inverse iteration}} (MP SII).
The approach is first illustrated on computation of a single eigenvector (Section \ref{sec:iii1}) and then extended to the block case (Section \ref{sec:iib}).

\item we calculate vibrational spectra of acetonitrile molecule CH$_3$CN using the proposed approach (Section \ref{sec:12}). The results are {more accurate than those of  the Smolyak grid approach \cite{ac-smolvibr-2011} but with much less storage requirements, and more accurate than the recent memory-efficient H-RRBPM method \cite{tc-hrrbpm-2015}, which is also based on tensor decompositions.
We note that the Smolyak grid approach does not require PES to be approximated by a small number of sum-of-product functions.}

\end{itemize}

\section{Discretization} \label{sec:discr}

We follow \cite{tc-hrrbpm-2015} and consider Schr\"odinger equation with omitted $\pi - \pi$ cross terms and the potential-like term in the normal coordinate kinetic energy operator. The Hamiltonian in this case looks as
\begin{equation}\label{eq:schroed}
	\mH = -\frac{1}{2}\sum_{i=1}^d \omega_i \frac{\partial^2}{\partial q_i^2}+ V(q_1,\dots, q_d),
\end{equation}
where $V$ denotes potential energy surface (PES).

We discretize the problem \eqref{eq:schroed} using the discrete variable representation (DVR) scheme on the tensor product of Hermite meshes
\cite{baye-dvr-1986} such that each unknown eigenfunction is represented as
\begin{equation}\label{eq:vibr_discr}
	\Psi_k(q_1,\dots,q_d) \approx \sum_{i_1, \dots, i_d=1}^n \mathcal{X}_{i_1,\dots,i_d}^{(k)}\, \psi_{i_1}(q_1) \dots \psi_{i_d}(q_d),
\end{equation}
where $\psi_{i}(q_i)$ denotes one-dimensional DVR basis function.
We call arising multidimensional arrays \emph{tensors}.

The Hamiltonian \eqref{eq:schroed} consists of two terms: the Laplacian-like part and the PES.
It is well-known that the Laplacian-like term can be written in the Kronecker product form
$$
\mathcal{D} =  D_1\otimes I\otimes\dots\otimes I + \dots + I\otimes\dots\otimes I\otimes D_d,
$$
where $D_i$ is the one dimensional discretization of the $i$-th mode.

The DVR discretization of the PES is represented as a tensor $\mathcal{V}$.
The operator corresponding to the multiplication by $\mathcal{V}$ is diagonal.
Finally the Hamiltonian is written as
\begin{equation}\label{eq:hamiltonian}
	\mathcal{H} = D_1\otimes I\otimes\dots\otimes I + \dots + I\otimes\dots\otimes I\otimes D_d + \text{diag}\left(\mathcal{V}\right).
\end{equation}

For our purposes it is convenient to treat~$\mathcal{H}$ not as a 2D $ \mathbb{R}^{n^d} \to  \mathbb{R}^{n^d}$ matrix, but as a multidimensional operator $\mathcal{H}: \mathbb{R}^{n\times\dots \times n} \to  \mathbb{R}^{n\times\dots \times n}$.
In this case the discretized Schr\"odinger equation has the form
\begin{equation}\label{eq:discr}
\sum_{j_1,\dots, j_d=1}^n \mathcal{H}_{i_1,\dots,i_d, j_1,\dots, j_d} \mathcal{X}_ {j_1,\dots, j_d} = E \mathcal{X}_ {i_1,\dots, i_d}.
\end{equation}
Hereinafter we use notation $\mH (\mX)$ implying matrix-by-vector product from \eqref{eq:discr}. Using this notation \eqref{eq:discr} can be equivalently written as $$\mH (\mX) = E \mX.$$

\section{Computing a single eigenvector}\label{sec:sing}

In this section we discuss how to solve the Schr\"odinger equation \eqref{eq:discr} numerically and  present our approach for finding a single eigenvector.
The case of multiple eigenvalues is discussed in Section \ref{sec:block}.

The standard way to find required number of smallest eigenvalues is to use iterative methods.
The simplest iterative method of finding one smallest eigenvalue is the shifted power iteration
\[
\begin{split}
	&\mX_{k+1} = (\mH - \sigma \mathcal{I}) \mX_k, \\
	&\mX_{k+1} := \mX_{k+1} / \sqrt{\left<\mX_{k+1}, \mX_{k+1} \right>},
\end{split}
\]
where the shift $\sigma$ is an approximation to the largest eigenvalue of~$\mH$.
The matrix-by-vector product is the bottleneck operation in this procedure.
This method was successfully applied to the calculation of vibrational spectra in \cite{lc-rrbpm-2014}. The eigenvectors in this work are represented as sum-of-products, which allows for fast matrix-by-vector multiplications.
Despite the ease of implementation and efficiency of each iteration, the convergence of this method requires thousands of iterations.

Instead of power iteration we use inverse iteration~\cite{ipsen-inverse-1997}
\begin{equation}\label{eq:inverse}
\begin{split}
	&\mX_{k+1} = (\mH - \sigma \mathcal{I})^{-1} \mX_k, \\
	&\mX_{k+1} := \mX_{k+1} / \sqrt{\left<\mX_{k+1}, \mX_{k+1} \right>},
\end{split}
\end{equation}
which is known to have much faster convergence if a good approximation $\sigma$ to the required eigenvalue $E^{(1)}$ is known.
{Question related the solution of linear systems in the TT-format will be discussed in Section \ref{sec:iii1}.}
Convergence of the inverse iteration is determined by ratio
\begin{equation}\label{eq:rho}
	\rho = \left| \frac{E^{(1)} - \sigma}{E^{(2)} - \sigma} \right|,
\end{equation}
where $E^{(2)}$ is the next closest to $E^{(1)}$ eigenvalue.
Thus, $\sigma$ has to be closer to $E^{(1)}$ than to $E^{(2)}$ for the method to converge. However, the closer $\sigma$ to $E^{(1)}$ is, the more difficult to solve the linear system with matrix $(\mathcal{H} - \sigma \mathcal{I})$ it is.
Therefore, typically this system is solved inexactly~\cite{golub-inexact-2000,berns-inexact-2006}.
Parameter $\sigma$ can also depend on the iteration number (Rayleigh quotient iteration), however in our experiments (Section \ref{sec:num}) constant choice of $\sigma$ yields convergence in 5 iterations.

As it follows from \eqref{eq:rho} for the inverse iteration to converge fast a good initial approximation has to be found.
To get initial approximation we propose to use locally optimal block preconditioned conjugate gradient (LOBPCG) method as it is relatively easy to implement in tensor formats, and a preconditioner can be explicitly utilized.

\subsection{TT-format}\label{sec:tt}
The problem with straightforward usage of the iterative processes above is that we need to store an eigenvector~$\mX_k$.
The storage of this vector is $\mO(n^d)$, which is prohibitive even for $n=15$ and $d=10$.
Therefore we need a compact representation of an eigenfunction which allows to do inversion and basic linear algebra operations in a fast way.
For this purpose we use the tensor train (TT, MPS) format \cite{osel-tt-2011}.
Tensor $\mX$ is said to be in the TT-format if it is represented as
\begin{equation}\label{eq:tt}
	\mathcal{X}_{i_1,i_2,\dots,i_d} \triangleq X^{(1)}(i_1)\, X^{(2)}(i_2)\, \dots\, X^{(d)}(i_d),
\end{equation}
where $X^{(k)}(i_k)  \in \mathbb{R}^{r_{k-1}\times r_k}$, $r_0 = r_d = 1$, $i_k = 1,\dots,n$.
Matrices $X_k(i_k)$ are called TT-cores and $r_k$ are called TT-ranks.
For simplicity we assume that $r_k=r$, $k=2,\dots,d-1$ and call $r$ the TT-rank.
{In numerical experiments we use different mode sizes $n_i$, $i=1,\dots,d$, but for simplicity we will use notation $n=\max_i n_i$ in complexity estimates.}
Compared to $\mO(n^d)$ parameters of the whole tensor, TT decomposition \eqref{eq:tt} contains $\mathcal{O}(dnr^2)$ parameters as each $X^{(i)}$, $i=2,\dots,d-1$ has size $n\times r \times r$.

The definition of the TT-format of an operator is similar to the TT representation of tensors
$$
\mathcal{H}_{i_1,\dots,i_d,\, j_1,\dots,j_d} \triangleq H^{(1)}(i_1,j_1)\, \dots\, H^{(d)}(i_d,j_d),
$$
where TT-cores $X^{(k)}(i_k, j_k)  \in \mathbb{R}^{R_{k-1}\times R_k}$, $R_0 = R_d = 1$.
If $R_k=r$, $k=2,\dots,d-1$ then this representation contains $\mathcal{O}(dn^2 R^2)$ degrees of freedom.

{
TT-format can be considered as a multidimensional generalization of SVD.
Other alternatives to generalize the SVD to many dimensions are the canonical, Tucker and Hierarchical Tucker formats. We refer the reader to \cite{kolda-review-2009,larskres-survey-2013} for detailed survey.
The important point why we use the TT decomposition is that it can be computed in a stable way and it does not suffer from the ``curse of dimensionality''.
Moreover, there exists well-developed software packages to work with the TT-decomposition~\cite{tt-toolbox}.
}

\subsection{Manifold-projected inverse iteration}\label{sec:iii1}

For the inverse iteration \eqref{eq:inverse} we need to find TT-representation of $\mX_{k+1}$ by approximately solving a linear system
\begin{equation}\label{sec:iieq}
(\mH - \sigma \mathcal{I}) \mX_{k+1} \approx \mX_k.
\end{equation}
Assume that both the exact eigenvector $\mX$ and the current approximation $\mX_k$  belong to the manifold $\mathcal{M}_r$ of tensors with TT-rank $r$.
The solution of \eqref{sec:iieq} may have ranks larger than $r$ and therefore be out of $\mathcal{M}_r$.
In the present work we suggest exploiting the information that $\mX$ belongs to $\mathcal{M}_r$ and retract the solution of \eqref{sec:iieq} back to the manifold.
We refer to this concept as a \emph{manifold-projected inverse iteration} (MP II).

In this work we pose the following optimization problem with a rank constraint
\begin{equation}\label{eq:min}
\begin{split}
&\underset{\mY}{\text{minimize}} \quad \ \mathcal{J}(\mY) \equiv \|(\mH - \sigma \mathcal{I}) \mathcal{Y} - \mX_k \|,\\
&\text{subject to}\quad \text{TT-rank}(\mY) = r.
\end{split}
\end{equation}
Problem \eqref{eq:min} is hard to solve as operator $(\mH - \sigma \mathcal{I})$ is close to singular.
Similarly to the inexact inverse iteration framework we are not searching for the solution that finds global minima of  \eqref{eq:min}, but utilize several sweeps of the alternating least squares (ALS) method with the initial guess $\mX_k$ \cite{holtz-ALS-DMRG-2012}.
The ALS procedure alternately fixes all but one TT-core and solves minimization problem with respect to this TT-core.
For instance, an update of a core $X^{(m)}$  when all other cores are fixed is $X^{(m)}_\text{new}$ found from
\begin{equation}\label{eq:alsmin}
X^{(m)}_\text{new} = \underset{X^{(m)}}{\text{arg\, min}}\ \mathcal{J}(X^{(1)},\dots, X^{(m)}, \dots, X^{(d)}).
\end{equation}
The minimization over a single core {(see Appendix \ref{sec:als})} is a standard linear least squares problem with the unknown vector that has the size $nr^2$ -- size of the corresponding TT-core and is very cheap. Moreover these systems can be also solved iteratively.
The described minimization over all cores in the TT-representation is referred to as one sweep of the ALS. {The total computational cost of one sweep is then~$\mathcal{O}(dn^2r^2 R^2)$, where $R$ is maximum TT-rank of the operator $\mathcal{H}$, see Appendix \ref{sec:als}.}

\emph{According to the proposed concept we start from $\mX_{k}$ and use only a few sweeps (one or two) of the ALS method  rather than running the method till convergence. Moreover, we found that one can solve local linear systems inexactly either with fixed number of iterations or fixed low accuracy.} 

Such low requirements for solution of local linear systems and number of ALS sweeps results in a very efficient method.
However, for this method to converge, a good initial approximation to both eigenvector and eigenvalue has to be found.

\subsection{Manifold-preconditioned LOBPCG for a single eigenvector}\label{sec:mp1}

To get initial approximation we use LOBPCG method.
The LOBPCG algorithm for one eigenvalue looks as follows
\begin{equation}\label{eq:lob}
\begin{split}
	& \mR_k := \mathcal{B} (\mH (\mX_k) - E_k \mX_k), \\
	&\mP_{k+1} :=  \alpha_2  \mR_k + \alpha_3 \mP_k, \\
	& \mX_{k+1} := \alpha_1 \mX_{k} + \mP_{k+1} ,  \\
	&\mX_{k+1} := \mX_{k+1} / \sqrt{\left<\mX_{k+1}, \mX_{k+1} \right>},
\end{split}
\end{equation}
where $\mathcal{B}$ denotes preconditioner and vector of coefficients ${ \alpha} = (\alpha_1, \alpha_2, \alpha_3)^T$ is chosen from minimization of the Rayleigh quotient
$$
R(\mX_{k+1}) = \frac{\left(\mH (\mX_{k+1}), \mX_{k+1}\right)}{\left(\mX_{k+1}, \mX_{k+1}\right)}.
$$
 Finding $\alpha$ is equivalent to solving the following $3\times 3$ eigenvalue problem
\begin{equation}\label{eq:rayleigh}
	 \begin{bmatrix} \mX_{k} \\ \mR_k \\ \mP_k \end{bmatrix}
	\mathcal{H} [\mX_{k}, \mR_k, \mP_k]\,
	\alpha = \lambda \begin{bmatrix} \mX_{k} \\ \mR_k \\ \mP_k \end{bmatrix} [\mX_{k}, \mR_k, \mP_k]\, \alpha.
\end{equation}

Let us discuss TT version of the LOBPCG.
Operations required to implement the LOBPCG are presented below:
\paragraph{Preconditioning.} The key part of the LOBPCG iteration is multiplication by a preconditioner $\mathcal{B}$.
{
In this work we use $\mathcal{B}\approx(\mathcal{H} - \sigma \mathcal{I})^{-1}$ as a preconditioner.
This preconditioner works well if the density of states is low, see~\cite{pc-chemprec-2001}.
To make a preconditioner more efficient one can project it to the orthogonal complement of current approximation of the solution, see Jacobi-Davidson method~\cite{svh-jd-2000}.
}

Instead of forming $(\mathcal{H} - \sigma \mathcal{I})^{-1}$ we calculate matrix-by-vector multiplication {$\mY \approx\mathcal{B} (\mH (\mX_k) - E_k \mX_k)$}. Hence, similarly to the inverse iteration (Section \ref{sec:iii1}) we propose solving a minimization problem
\[
\begin{split}
&\underset{\mY}{\text{minimize}} \quad\ \,  \|(\mathcal{H} - \sigma \mathcal{I})\mathcal{Y} -   (\mH (\mX_k) - E_k \mX_k) \|,\\
&\text{subject to}\quad \text{TT-rank}(\mY) = r.
\end{split}
\]
We also use only several sweeps of ALS for this problem.
We refer to this construction of preconditioner as a \emph{manifold preconditioner} as it retracts the residual on a manifold of tensors with fixed rank $r$.
Note that if $(\mathcal{H} - \sigma \mathcal{I})$ is known to be positive definite, then minimization of energy functional
{
$$
\mathcal{J}(\mY) = \left<(\mathcal{H} - \sigma \mathcal{I})\mathcal{Y}, \mY \right> - 2 \left< \mH (\mX_k) - E_k \mX_k, \mY\right>
$$}
can be used instead of minimization of the residual.

\paragraph{Summation of two tensors.}
Given two tensors $\mX$ and $\mY$ with ranks $r$ in the TT-format
\[
\begin{split}
	&\mX_{i_1\dots i_d} = X^{(1)}(i_1)\dots X^{(d)} (i_d), \\
	&\mY_{i_1\dots i_d} = Y^{(1)}(i_1)\dots Y^{(d)} (i_d),
\end{split}
\]
the cores of the sum $\mathcal{Z} = \mX + \mY$ are defined as \cite{osel-tt-2011}
$$
	Z^{(k)}(i_k) = \begin{bmatrix}  X^{(k)}(i_k) & 0 \\ 0 & Y^{(k)}(i_k) \end{bmatrix}, \quad k = \overline{2,d-1}
$$
and
$$
Z^{(1)}(i_1) = \begin{bmatrix} X^{(1)}(i_1) & Y^{(1)}(i_1) \end{bmatrix},\
Z^{(d)}(i_d) = \begin{bmatrix} X^{(d)}(i_d) \\ Y^{(d)}(i_d) \end{bmatrix},
$$
Thus, tensor $\mathcal{Z}$ is explicitly represented with ranks~$2r$.

\paragraph{Inner product and norm.}
{
To find inner product of two tensors $\mX$ and $\mY$ in the TT-format we first need to calculate the Hadamard product, which can calculated as
\[
\begin{split}
(\mX \circ \mY)_{i_1\dots i_d} = X^{(1)}(i_1)\dots X^{(d)} (i_d) Y^{(1)}(i_1)\dots Y^{(d)} (i_d)= \\  \left(X^{(1)}(i_1) \otimes Y^{(1)}(i_1)\right)\cdot \dots \cdot \left(X^{(d)}(i_d) \otimes Y^{(d)}(i_d)\right).
\end{split}
\]
Therefore,
$$
\left< \mX, \mY \right> = \Gamma_1 \dots \Gamma_d,
$$
where $$\Gamma_k = \sum_{i_k} X^{(k)}(i_k) \otimes Y^{(k)}(i_k).$$
Using special structure of matrices $\Gamma_k$ the inner product can be calculated in $\mathcal{O}(dnr^3)$ complexity \cite{osel-tt-2011}. The norm can be computed using inner product as $\|\mX\| = \sqrt{\left< \mX, \mX \right>}$.
}
\paragraph{Reducing rank (rounding).} As we have seen, after summation ranks grow. To avoid rank growth there exists special \emph{rounding} operation.
It suboptimally reduces rank of a tensor with a given accuracy.
{In 2D the rounding procedure of $\mX=UV^\top$ looks as follows. First we calculate QR decompositions of matrices $U$ and $V$: $$U=Q_U R_U, \quad V=Q_V R_V.$$ Hence
$$ \mX = Q_U R_U R_V^\top Q_U^\top.$$
Finally, to reduce the rank we calculate the SVD decomposition of $R_U R_V^\top$ and truncate singular values up to required accuracy.
In \cite{osel-tt-2011} this idea is generalized to the TT case. The complexity is $\mathcal{O}(dnr^3)$.
}

\begin{figure*}
\begin{minipage}{\linewidth}
\begin{algorithm}[H]
  \caption{Auxiliary functions}
  \label{alg:functions}
\begin{itemize}
\item[] $ \mY= \verb MULTIFUNCRS (\verb|func|, \bX) $: calculates $\verb|func|(\bX)$ via the cross approximation algorithm (Section \ref{sec:iib}).
\item[] $ \bY= \verb BLOCK_MATVEC (\bX, M)$: block multiplication of a vector of TT tensors $\bX$ by real-valued matrix $M$ using $\verb MULTIFUNCRS $ function, see \eqref{eq:bmatvec}.
\item[] $\bY = \verb QR (\bX)$: orthogonalizes TT tensors $\mX_1,\dots,\mX_B$: $\bX=(\mX_1,\dots,\mX_B)$ via Cholesky (Section~\ref{sec:iib}) or modified Gram-Schmidt procedure.
Matrix-by-vector multiplications are done using using $\verb BLOCK_MATVEC $.
\item[] $\bX = \verb ALS (\mA,\ {\bf F},\ n_\text{swp})$: solves $\mA (\mX^{(i)}) = \mathcal{F}^{(i)}$, $i=1,\dots,$\verb|length|$(\bF)$ using $n_\text{swp}$ sweeps of ALS method to minimize $\|\mA (\mX^{(i)}) - \mathcal{F}^{(i)}\|$ with a rank constraint.
\item[] $\bY = \verb ORTHO (\bX, \bQ)$: orthogonalizes TT tensors $\mX_1,\dots,\mX_B$: $\bX=(\mX_1,\dots,\mX_B)$ with respect to $\bQ$: $\mY^{(i)} = \mX^{(i)} - \sum_{j=1}^B \left<\mX^{(i)}, \mQ^{(j)} \right> \mQ^{(j)} $. To avoid rank growth we use rounding if \verb|length|$(\bX)$ is small and via \verb|MULTIFUNCRS| if \verb|length|$(\bX)$ is large.
\item[] $\bY = \Tau_r(\bX)$: truncates each tensor $\mX_1,\dots,\mX_B$: $\bX=(\mX_1,\dots,\mX_B)$ with rank~$r$ using rounding procedure.
\end{itemize}
\end{algorithm}
\end{minipage}
\end{figure*}

\paragraph{Matrix-by-vector multiplication.} To calculate matrix-by-vector operation $\mathcal{H}(\mX)$ in \eqref{eq:lob} it is convenient to have Hamiltonian $\mH$ represented in the TT-format.
In this case there exists explicit representation of matrix-by-vector multiplication $\mY = \mH (\mX)$ when both $\mathcal{H}$ and $\mX$ are represented in the TT-format \cite{osel-tt-2011}
$$
Y^{(k)}  ({i_k}) = \sum_{j_k} \left( H^{(k)}(i_k,j_k) \otimes X^{(k)}(j_k) \right),
$$
{which gives representation with TT-rank~$=rR$.}
To reduce the rank one can either use the rounding procedure or use ALS minimization of the following optimization problem
\[
\begin{split}
&\underset{\widehat\mY}{\text{minimize}}\quad\ \|\widehat\mY - \mY\|,\\
&\text{subject to}\quad \text{TT-rank}(\widehat\mY) = r.
\end{split}
\]
which is faster than rounding for large ranks.

\section{The block case} \label{sec:block}

In the previous section we discussed the algorithm for a single eigenvector.
In this section we extend the algorithm to a block case.
The difference is that we need to make additional block operations such as othogonalization and block multiplication.

\subsection{Manifold-projected simultaneous inverse iteration}\label{sec:iib}

{In this section we present manifold-projected version of the simultaneous inverse iteration, which is also known as block inverse power method.}

Assume we are given a good initial approximation to eigenvalues and eigenvectors of linear operator $\mathcal{H}$ given by its TT-representation.
In this case inverse iteration yields fast convergence rate with low memory requirements.

\begin{algorithm}[H]
  \caption{Manifold-Projected Simultaneous Inverse Iteration (MP SII) with shifts}\label{alg:inverse}
\label{alg:lobpcg}
  \begin{algorithmic}[1]
\Require{TT-matrix $\mathcal{H}$; initial guess $\bX_{0} = [\mathcal{X}_0^{(1)} \dots \mathcal{X}_0^{(B)}]$ where $\mX_0^{(i)}$ are TT tensors; truncation rank~$r$}
\Ensure{$\Lambda$ and $\bX$ -- approximation of $B$ eigenvalues that are close to $\sigma$ and corresponding eigenvectors}
\Statex
\State Group close eigenvalues in clusters
\State $X := [\,]$, $\Lambda := [\,]$
\For{each cluster with index $\nu$}
\State Compute shift $\sigma_\nu$ -- average eigenvalue in cluster $\nu$
\State Compute $\bX_{0}^\nu$ -- corresponding subvector of $\bX_{0}$
\State{$\bX_{0}^\nu :=  \verb QR (\bX_{0}^\nu)$}
	\For{$k=0,1,\dots$ until converged}
		\State $\bX_{k+1}^\nu := \verb ALS (\mathcal{H} - \sigma_\nu \mathcal{I},\ \bX_k^\nu,\ n_\text{swp}$)
		\State{$\bX_{k+1}^\nu :=  \verb QR (\bX_{k+1}^\nu)$}
	\EndFor
\State $\bX: = [\bX, \bX_{k+1}^\nu]$, $\Lambda: = [\Lambda, \text{diag}({\bX_{k+1}^\nu}^T \mA (\bX_{k+1}^\nu)]$
\EndFor
\Statex
\Return $\Lambda,\bX$
  \end{algorithmic}
\end{algorithm}

Given initial approximation we first split eigenvalues into clusters of eigenvalues. Proximity of eigenvalues is defined using a threshold parameter.
If a cluster consists of one eigenvalue, then we run a version described in Section~\ref{sec:iii1}.
Otherwise we need additional orthogonalization on each iteration step.
The orthogonalization is done using the QR decomposition via Cholesky factorization.
Let us consider orthogonalization procedure for vector $\bX=[\mX^{(1)}, \dots, \mX^{(B)}]$ in more details.
First we calculate Gram matrix
$$
G = \bX^T \bX,
$$
where we used notation
$$
	\bX^T \bY \triangleq
\begin{bmatrix}
	\left< \mX^{(1)}, \mY^{(1)}\right> & \dots & \left< \mX^{(1)}, \mY^{(B)}\right> \\
	\vdots & & \vdots\\
	\left< \mX^{(B)}, \mY^{(1)}\right> & \dots & \left< \mX^{(B)}, \mY^{(B)}\right> \\
\end{bmatrix}.
$$
The calculation of Gram matrix can be done in $\mO(B^2ndr^3)$ operations {as calculation of a single inner product requires $\mO(ndr^3)$ operations.}
Then we calculate the Cholesky factorization of the $B\times B$ Gram matrix $G = LL^T$.
{We note that $G$ consists of numbers, so the standard Cholesky factorization of matrices is used.}
The final and the most time consuming step is to find block matrix-by-vector operation denoted by \verb|BLOCK_MATVEC|$(\bX, L)$. For a general matrix $M\in\mathbb{R}^{P\times B}$ function \verb|BLOCK_MATVEC|$(\bX, M)$ produces block vector $\bY=[\mY^{(1)}, \dots, \mY^{(P)}]$ such that
\begin{equation}\label{eq:bmatvec}
\mY^{(i)} = \sum_{j=1}^P L_{ij} {\mX}^{(j)},\quad i=1,\dots,B.
\end{equation}
If $P$ is small, say $P<20$, then the summation can be done using summation and rounding.
{The typical value of $P$ in numerical experiments is $B=80$ or $2B=160$, so to reduce the constant in complexity we use the so-called} \emph{cross approximation method}, which is able to calculate the TT-decomposition of a tensor using only few of its elements.
Namely, the cross approximation is able to build TT-representation {of a tensor of exact rank~$r$} using only $\mO(dnr^2)$ elements (number of parameters in the TT-format) using interpolation formula \cite{ot-ttcross-2010} with $\mO(dnr^3)$ complexity.
{If tensor is approximately of low rank then tensor can be approximated with given accuracy  using the same approach and quasioptimality estimates exist \cite{sav-qott-2014}.}
To find TT-representation of $\mY^{(i)}$ we calculate $\mO(dnr^2)$ elements of tensor $\mY^{(i)}$ by explicitly calculating elements in ${\mX}^{(j)}$ and summing them up with coefficients $L_{ij}$.
This approach allows to calculate general type of functions of a block vector $f(\bX)$ and is referred to as \verb|MULTIFUNCRS|. {It is typically used when either the number of input tensors is large or some nonlinear function of input tensors must be computed.}
{The similar approach was used in \cite{ro-crossconv-2015} for a fast computation of convolution.}

The result of solving the block system $\mH(\bX) = \bF$ using $n_\text{swp}$ of optimization procedure is denoted as $\bX = \verb|ALS|(\mH,\ {\bf F},\ n_\text{swp})$, where $$\mH (\bX) \triangleq [\mH(\mX^{(1)}), \dots, \mH (\mX^{(B)})].$$  

The overall algorithm is summarized in Algorithm \ref{alg:inverse}. All discussed auxiliary functions are presented in Algorithm \ref{alg:functions}.
If the cluster size is much smaller than the number of required eigenvalues, the complexity of finding each cluster is fully defined by the inversion operation which costs $\mO(dn^2r^2 R^2)$.
Thus, the overall complexity of the inverse iteration scales linearly with $B$: $\mO(Bdn^2r^2 R^2)$.
We note that each cluster can be processed independently in parallel.


\subsection{Manifold-preconditioned LOBPCG method}\label{sec:mpb}

We use the LOBPCG method to get initial approximation for the manifold-preconditioned simultaneous inverse iteration.
The problem is that each iteration of the LOBPCG is much slower compared to the inverse iteration for large number of eigenvalues.
We hence only run LOBPCG with small ranks and then correct its solution with larger ranks using the inverse iteration.

\begin{algorithm}[H]\label{alg:lobpcg}
  \caption{Manifold-Preconditioned LOBPCG (MP LOBPCG)}\label{alg-lobpcg}
\label{alg:lobpcg}
  \begin{algorithmic}[1]
\Require{TT-matrix $\mathcal{H}$; initial guess $\bX_{0} = [\mathcal{X}_0^{(1)} \dots \mathcal{X}_0^{(B)}]$, where $\mX_0^{(i)}$ are TT tensors; truncation rank $r$}
\Ensure{$\Lambda$ and $\bX$ -- approximation of $B$ smallest eigenvalues and eigenvectors of $\mathcal{H}$}
\Statex
	\State{$\bX_{0} :=  \verb QR (\bX_{0})$}
	\State{$\left(\bX_{0}^T \mathcal{H} (\bX_{0})\right) {S}_0= {S}_0        \Lambda_{0}$\Comment{Eigendecomposition}}
	\State{$\bX_0 := \verb BLOCK_MATVEC (\bX_0, S_0)$}
	\State $\bR_0 :=\mathcal{H}(\bX_0) - \bX_0\Lambda_0$
	\State{$\bP_0 := 0\cdot \bX_0$}
	\For{$k=0,1,\dots$ until converged}
		\State $\bR_k := \verb ORTHO (\bR_k, \bQ)$  \Comment{$\bQ$ -- converged vectors}
		\State{$\bR_k := \verb ALS (\mathcal{H} - \sigma \mathcal{I},\ \bR_k,\ n_\text{swp}$)}			                 
		\State{$\widetilde H := [\bX_k, \bR_k, \bP_k]^T \mathcal{H}([\bX_k, \bR_k, \bP_k])$}
		\State{$\widetilde M := [\bX_k, \bR_k, \bP_k]^T [\bX_k, \bR_k, \bP_k]$}
		\State{$\widetilde H\widetilde {S}_k =\widetilde M \widetilde {S}_k \widetilde\Lambda_k, \widetilde {S}_k^T \widetilde M \widetilde {S}_k = I_{3B}$ \Comment{Eigendecomposition}}
		\State{$\bP_{k+1} := \verb BLOCK_MATVEC ([\bR_k, \bP_k],\ \widetilde {S}_k (B\text{:}3B,1\text{:}B) ) $}
		\State$\bX_{k+1} :=   \verb BLOCK_MATVEC ( \bX_k,\ \widetilde {S}_k(1\text{:}B,1\text{:}B)) +  \bP_{k+1}$
		\State  $\bX_{k+1}:= \Tau_r(\bX_{k+1})$
		\State{$\Lambda_{k+1} := \widetilde\Lambda_k (1\text{:}B,1\text{:}B)$}
		\State $\bR_{k+1} :=\mathcal{H}(\bX_{k+1}) - \bX_{k+1}\Lambda_{k+1}$
		\If{Some tensors from $\bR_{k+1}$ converged}
			\State Augment $\bQ$ with corresponding tensors from $\bX_{k+1}$
			\State Restart the algorithm with new $\bX_{0}$
			\State (Optionally) increase truncation rank
		\EndIf

	\EndFor
\Statex
\Return $\Lambda_{k+1},\bX_{k+1} $
  \end{algorithmic}
\end{algorithm}

The LOBPCG algorithm is summarized in Algorithm~\ref{alg:lobpcg}.
Auxiliary functions used in the algorithm are presented in Algorithm \ref{alg:functions}.
We also use MATLAB like notation for submatrices, e.g. $S(B:3B,1:B)$ is the corresponding submatrix in matrix $S$.

Block matrix-by-vector multiplication \eqref{eq:bmatvec} arises when multiplying by $B\times 2B$ matrix, where $B$ is the number of eigenvalues to be found.
When $2B$ is a large number and we use \verb|MULTIFUNCRS| for block matrix-by-vector product instead of straightforward truncation.
This is the most time consuming step in the algorithm and it costs $\mO(B^2 dnr^3)$.
Another time-consuming part is matrix-by-vector multiplication which costs $\mO(B dn^2 r^2 R^2)$. Thus, the overall complexity of each iteration is $\mO(B dnr^2 (Br +nR^2))$.

In order to accelerate the method, we consider a version of LOBPCG with deflation.
Deflation technique is that we stop iterating converged eigenvectors.
In this case the residual must be orthogonalized with respect to the converged eigenvectors.
This procedure is denoted as \verb|ORTHO| and is described in more details in Algorithm \ref{alg:functions}.
It might be also useful to increase rank of vectors that did not converged after previous deflation step.


\section{Numerical experiments} \label{sec:num}

The prototype is implemented in Python using the \verb|ttpy| library \url{https://github.com/oseledets/ttpy}.
The code of the proposed algorithm can be found at \url{https://bitbucket.org/rakhuba/ttvibr}.
For the basic linear algebra tasks the MKL library is used. Python and MKL are from the Anaconda Python Distribution
\url{https://www.continuum.io}. Python version is 2.7.11. MKL
version is 11.1-1. Tests were performed on a single Intel Core i7 2.6 GHz processor with 8GB of RAM.
{However, only 1 thread was used.}

\subsection{64-D bilinearly coupled oscillator}\label{sec:64}

First of all, we test our approach on a model Hamiltonian when the solution is known analytically.
Following~\cite{tc-hrrbpm-2015} we choose bilinearly coupled 64 dimensional Hamiltonian
$$
	\mathcal{H} = \sum_{i=1}^d \frac{\omega_i}{2} \left(-\frac{\partial^2}{\partial q_i^2} + q_i^2\right ) +  \sum_{j=1}^{d-1} \sum_{i>j} \alpha_{ij} q_i q_j,
$$
with $\omega_j = \sqrt{j/2}$ and $\alpha_{ij}=0.1$. {The TT-rank of this hamiltonian is $3$ for all inner modes independently of $d$ or $n$.}

\begin{figure}[t]
\includegraphics[width=9cm]{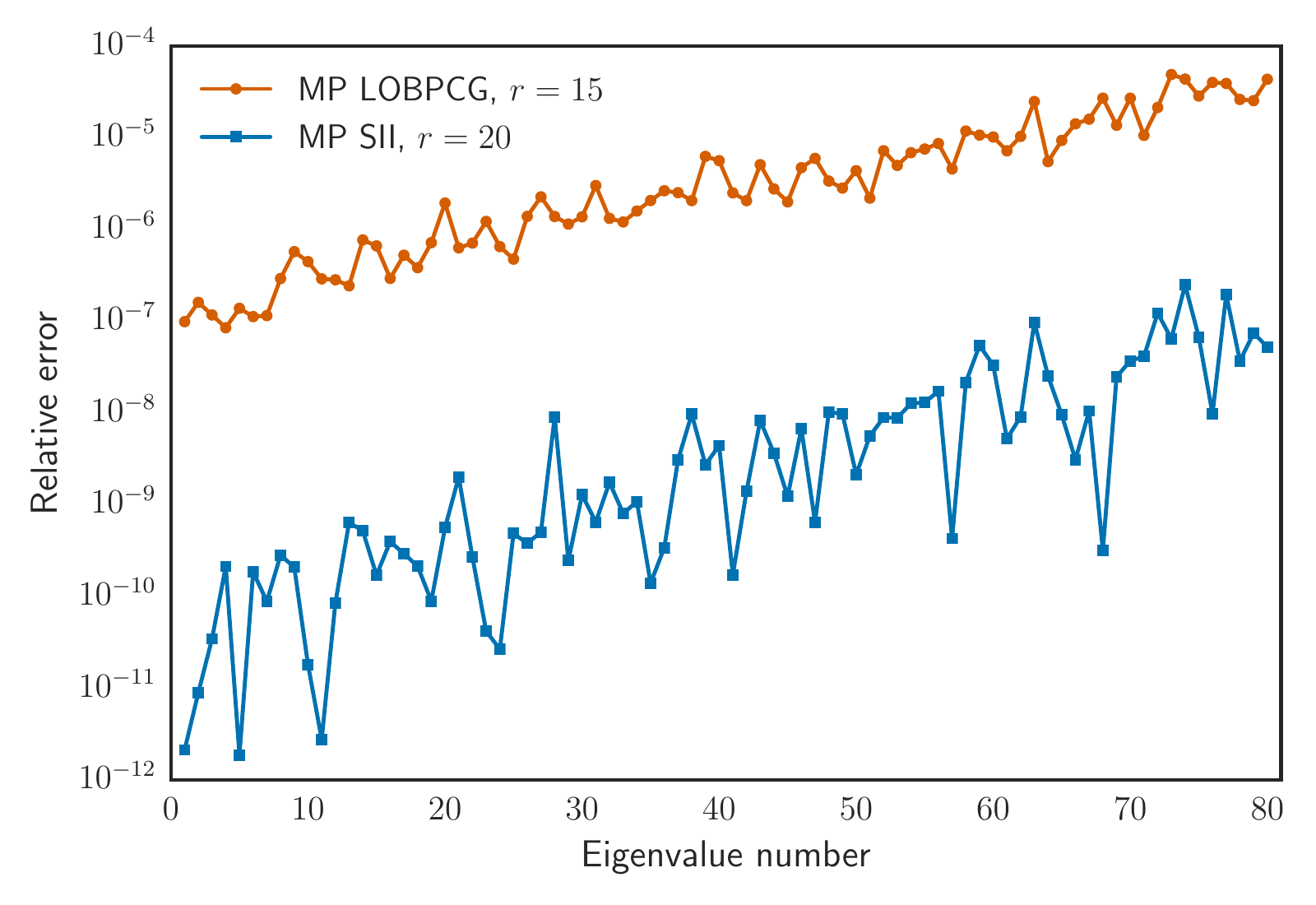}
\caption{Relative error with respect to the eigenvalue number for 64-D bilinearly coupled oscillator. MP LOBPCG stands for manifold-preconditioned LOBPCG and MP SII stands for manifold-projected simultaneous inverse iteration. MP SII uses solution of MP LOBPCG as an initial guess.} \label{fig:berror}
\end{figure}

To solve this problem we first use manifold-preconditioned LOBPCG method (Section~\ref{sec:mpb}) with rank $r=15$ and then correct it with the MP inverse iteration.
{The thresholding parameter $\delta$ for separation of energy levels into clusters in the inverse iteration is $10^{-4}$ ($E_i$ and $E_{i+1}$ are in the same cluster if $|E_i - E_{i+1}|<\delta \cdot |E_{i+1}|$).
The mode size $n=15$ is constant for each mode.}
As it follows from Figure \ref{fig:berror} the MP inverse iteration significantly improves the accuracy of the solution.
We used 10 iterations of the LOBPCG.
The MP LOBPCG computations took approximately 3 hours of CPU time and the MP inverse iteration took additionally 30 minutes.

\begin{figure}
\includegraphics[width=9cm]{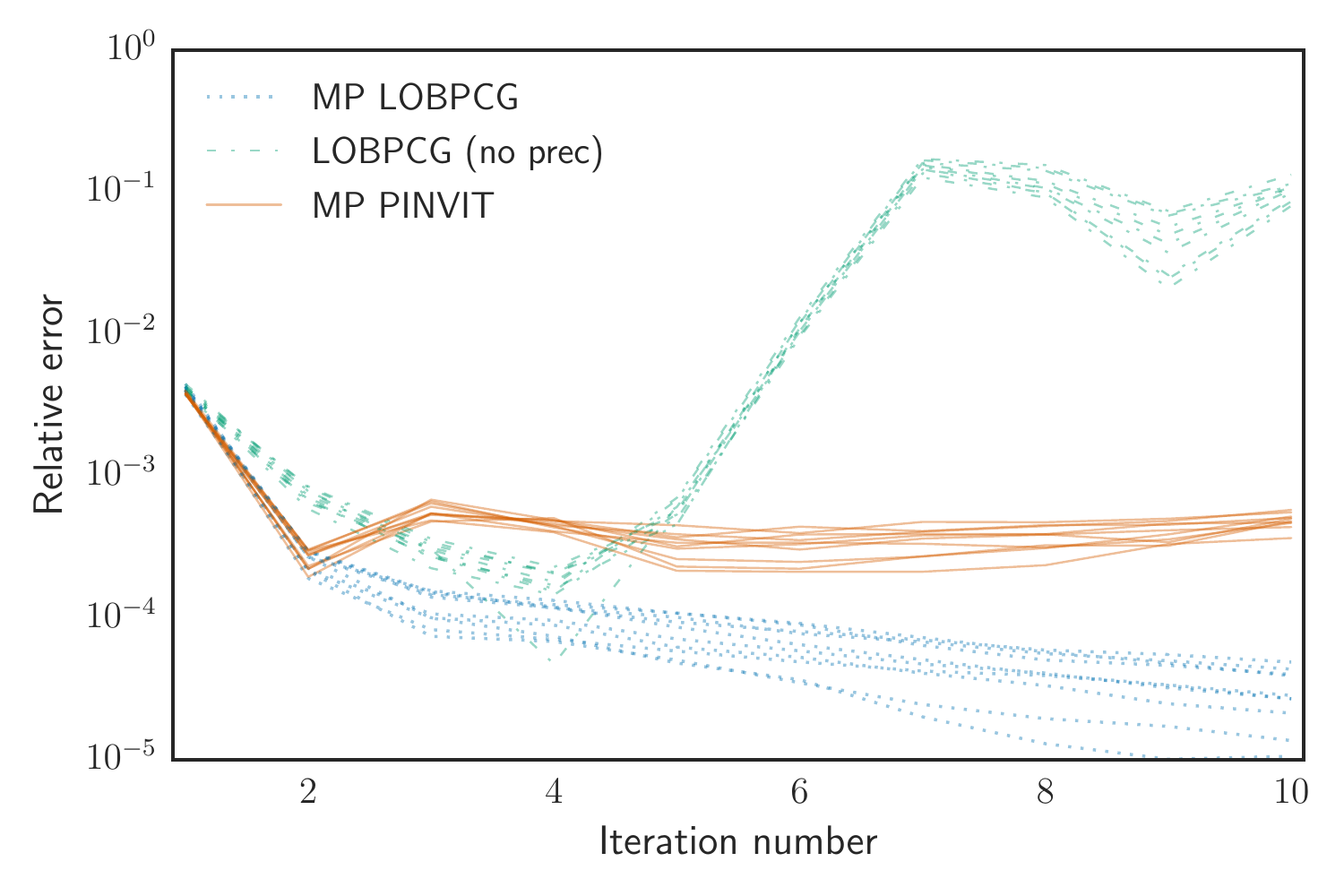}
\caption{
Relative error with respect to the iteration number for 64-D bilinearly coupled oscillator and different iterative methods. For each method the convergence of last 10 eigenvalues is presented. PINVIT denotes preconditioned inverse iteration \eqref{eq:pinvit}. MP stands for manifold preconditioner.} \label{fig:conv_sd}
\end{figure}

We also tested tensor version of the preconditioned inverse iteration (PINVIT) \cite{mach-innereig-2013} that in case of a single eigenvector looks as
\begin{equation}\label{eq:pinvit}
\begin{split}
	 &\mX_{k+1} = \Tau_r\left(\mX_k - \tau_k\, \mathcal{B} (\mH (\mX_k) -  E_k \mX_k)\right), \\
	&\mX_{k+1} := \mX_{k+1}/ \sqrt{\left<\mX_{k+1}, \mX_{k+1} \right>},
\end{split}
\end{equation}
where $\tau_k$ is selected to minimize the Rayleigh quotient.
Figure \ref{fig:conv_sd} illustrates convergence behavior of last 10 eigenvalues for different methods.
The PINVIT which also allows for explicit preconditioner converged to wrong eigenvalues. LOBPCG method without preconditioner is unstable due to rank thresholding.
We note that all these iterations converged to correct eigenvalues if the number of eigenvalues to be found was less than~$30$.

\subsection{Acetonitrile (CH$_3$CN) molecule}\label{sec:12}

In this part we present calculations of vibrational spectra of acetonitrile (CH$_3$CN) molecule.
The Hamiltonian used is described in  \cite{ac-smolvibr-2011} and looks as follows
\[
\begin{split}
V(q_1,\dots,q_{12}) = \frac{1}{2} \sum_{i=1}^{12} \omega_i q_i^2 + \frac{1}{6} \sum_{i=1}^{12} \sum_{j=1}^{12} \sum_{k=1}^{12} \phi_{ijk}^{(3)} q_i q_j q_k\\
+ \frac{1}{24} \sum_{i=1}^{12} \sum_{j=1}^{12} \sum_{k=1}^{12} \sum_{l=1}^{12} \phi_{ijkl}^{(4)} q_i q_j q_k q_l.
\end{split}
\]
It contains $323$ terms: $12$ kinetic energy terms, $12$ quadratic, $108$ cubic, and $191$ quartic potential terms.
{We chose the same basis size that was used in \cite{lc-rrbpm-2014}, namely the mode sizes were $\{9, 7, 9, 9, 9, 9, 7, 7, 9, 9, 27, 27 \}$ corresponding to the order described in that work.
We found that ranks of the Hamiltonian for this particular molecule do not strongly depend on the permutation of indices, namely the largest rank we observed among random permutations was $31$, while the maximum rank of the best permutation was $23$.
In computations we permuted indices such that array of $\omega_i$ is sorted in a decaying order.
Table \ref{tab:hamranks} contains ranks of the Hamiltonian in the TT-representation.
We note that  total ranks are ranks of a sum of potentials after rounding, and hence they are not equal to the sum of ranks of potentials in Table \ref{tab:hamranks}.
}

{
To assemble the potential $V$ one needs to add   rank-1 terms $q_i q_j q_k$ and rank-1 terms $q_i q_j q_k q_l$.
Each rank-1 term can be expressed analytically in the TT-format.
As was mentioned in Section \ref{sec:mp1} after each summation the rank grows, so the rounding procedure is used.
Recall that the rounding procedure requires $\mathcal{O}(d n^2 r^3)$ operations.
Thus, the complexity of assembling the Hamiltonian is
$$\mathcal{O}\left(\left(\verb nnz \left(\phi_{ijk}^{(3)}\right) + \verb nnz \left(\phi_{ijkl}^{(4)}\right) \right) d n^2 R^3 \right),$$
where \verb|nnz| stands for number of nonzeros, $n$ is the maximum mode size and $R$ is the maximum rank of $\mH$. The total time of assembling the Hamiltonian was less than $1$ second.
}

\setlength{\tabcolsep}{4pt}
\begin{table}[t!]
\caption{TT-ranks of the parts of the Hamiltonian with $\epsilon=10^{-10}$ threshold. \\}
\label{tab:hamranks}
\begin{tabular}{lcccccccccccc}
\hline
& $r_1$ & $r_2$ & $r_3$ & $r_4$ & $r_5$ & $r_6$ & $r_7$ & $r_8$ & $r_9$ & $r_{10}$ & $r_{11}$ \\
\hline
Quadratic & 2 & 2 & 2 &2 &2 &2 &2 &2 &2 &2 &2 \\
Cubic &  3 & 6& 11& 14& 14& 14& 14& 12& 9& 5& 3\\
Quartic &  5 & 7& 12& 19& 23& 26& 24& 18& 15& 8& 5\\
Total &  5& 9& 14&21&25&26&24&18&15&8&5\\
\hline
\end{tabular}
\end{table}

\begin{figure}[t]
\includegraphics[width=9cm]{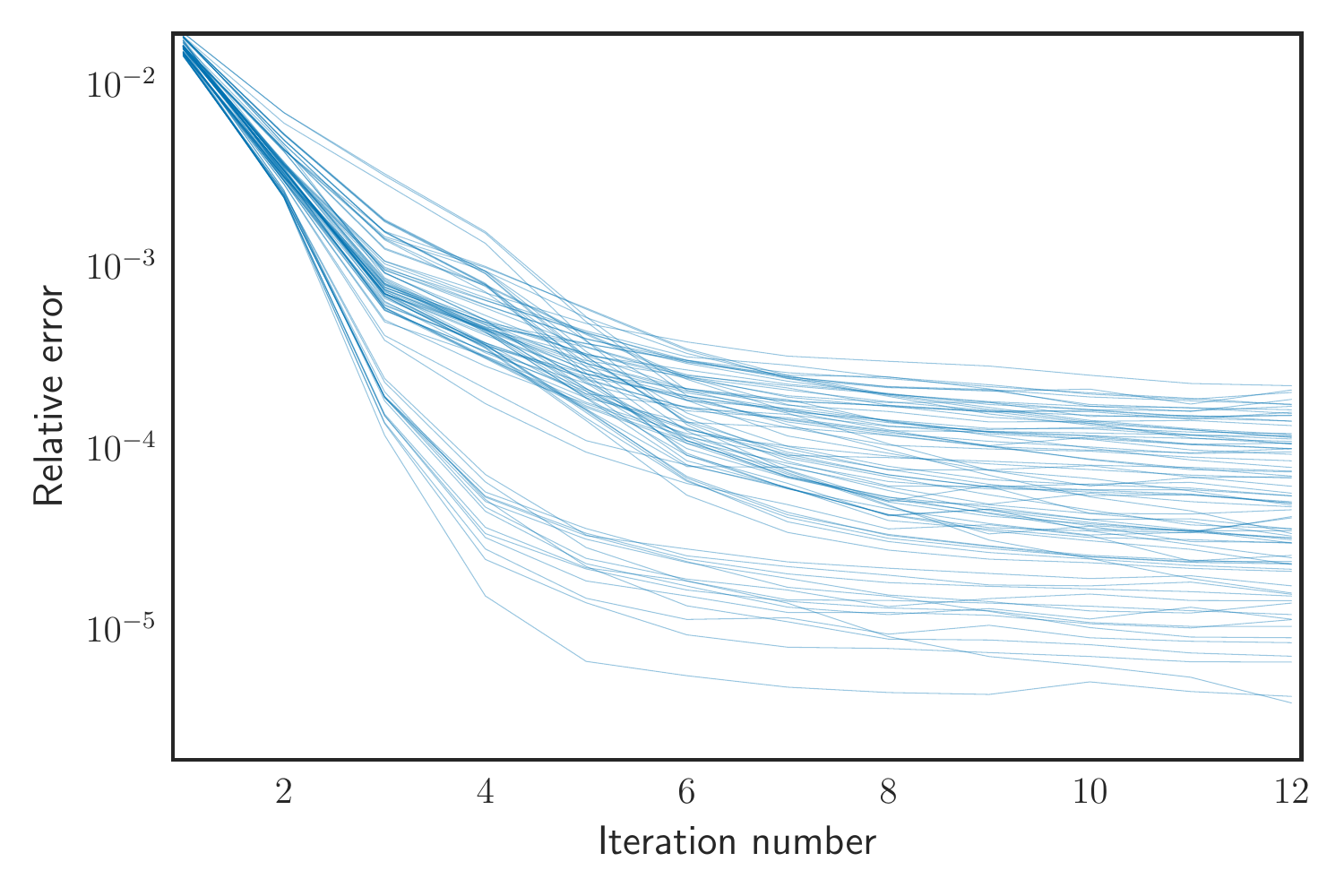}
\caption{Relative error of each of 84 eigenvalues of the acetonitrile molecule with respect to iteration number for the manifold-preconditioned LOBPCG. Relative error is calculated using Smolyak quadratures~\cite{ac-smolvibr-2011} as a reference value.} \label{fig:conv_ch3cn}
\end{figure}

\begin{figure}[t]
\includegraphics[width=9cm]{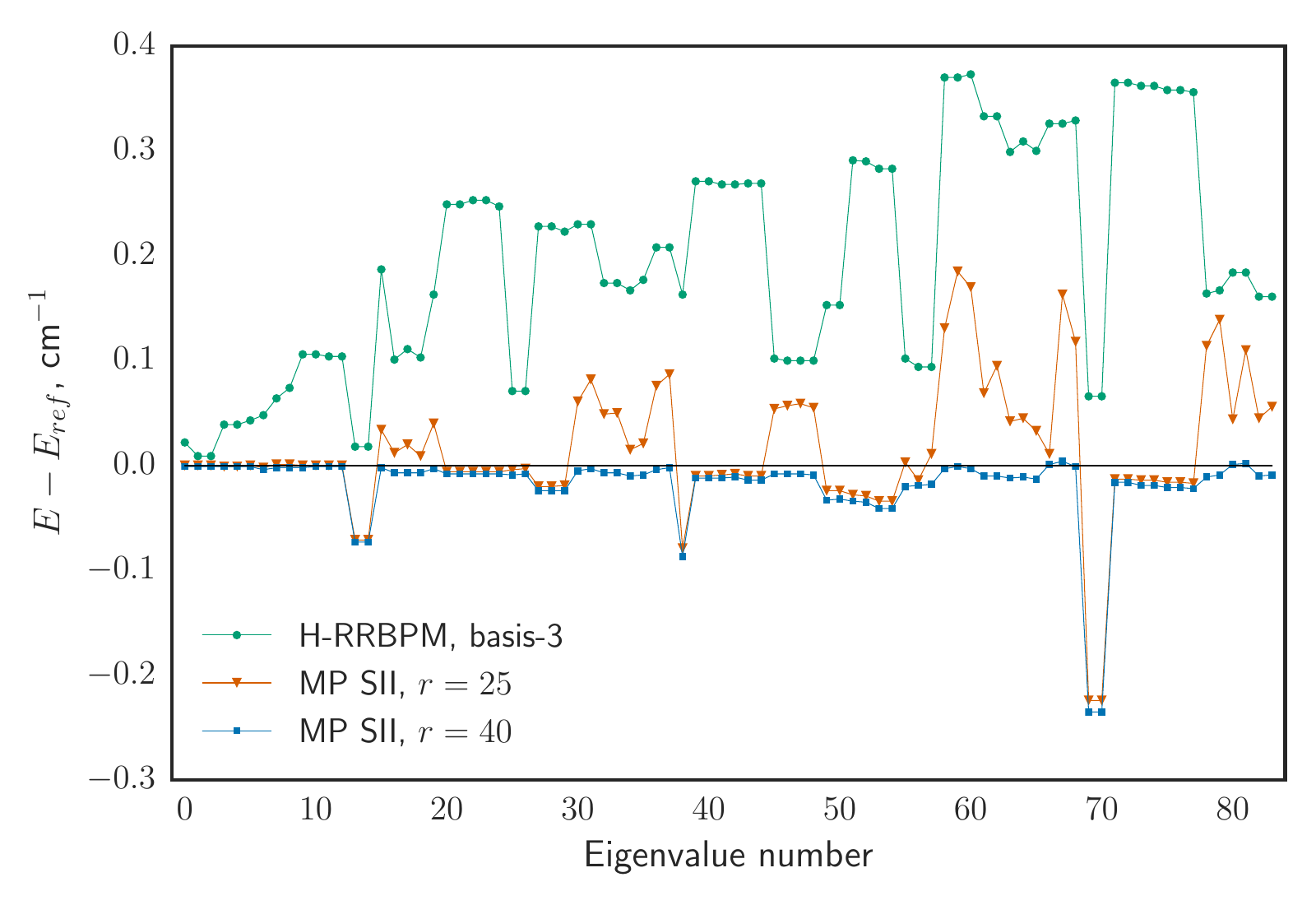}
\caption{Relative error $E - E_\text{ref}$ of eigenvalues of the acetonitrile molecule with respect to the eigenvalue number. The $E_\text{ref}$ are energies obtained by Smolyak quadratures \cite{ac-smolvibr-2011}. Negative value of error stand for more accurate than Smolyak quadratures energies. Black line denotes zero value of the error. MP SII stands for manifold-projected inverse iteration.} \label{fig:ch3cn_error}
\end{figure}

We ran the LOBPCG method with TT-rank equal to 12 and used manifold preconditioner.
Initial guess is chosen from the solution of the harmonic part of the Hamiltonian.
Eigenvectors of multidimensional quantum harmonic oscillator are tensor product of 1D oscillator eigenvectors and therefore can be represented analytically in the TT-format with rank $1$.
Shift for the preconditioner is set to be the lowest energy of the harmonic part.
The convergence of each eigenvalue is presented on Figure \ref{fig:conv_ch3cn}.
The obtained eigenfunctions were used as an initial approximation to the inverse iterations with ranks equal to~$25$.
Shifts were chosen to be LOBPCG energies.
{The thresholding parameter $\delta$ for separation of energy levels into clusters in the inverse iteration is $10^{-3}$}
These results were further corrected with the inverse iteration with rank $r=40$.
As it follows from Table \ref{tab:energies} and Figure~\ref{fig:ch3cn_error} both corrections are more accurate than the H-RRBPM method.
The latter correction with rank $r=40$ yields energy levels lower than those of the Smolyak quadrature method \cite{ac-smolvibr-2011} which means that our energy levels are more accurate. We note that storage of a solution with $r=40$ is less than the storage of the Smolyak method (180~MB compared with 1.5~GB).


Timings and storage of the H-RRBPM method were taken from  \cite{tc-hrrbpm-2015}.
On this example the state-of-the-art method \verb|eigb| \cite{dkos-eigb-2014} method converges within approximately several days.
The problem is that all basis functions are considered in one basis, which leads to large ranks.
Nevertheless, the \verb|eigb| becomes very efficient and accurate when small number of eigenpairs are required.

\section{Related work} \label{sec:related}

The simplest basis set for representing unknown eigenfunctions is the direct product (DP) of one-dimensional basis functions.
If a fast matrix-by-vector operation is given, then Krylov iterative methods are available and the only problem is the exponential storage requirements~$n^d$.
Alternatively one can prune a direct product basis \cite{at-pruned-2011,ben-pruned-2004,dt-pruned-2005} or use a basis that is a product of functions with more than one variable \cite{bb-contr-1991,bl-floppy-1989}.

In this work, we focus on DP basis and further reduce~$n^d$ storage by approximating unknown eigenvectors in the TT-format.
We refer the reader to \cite{kolda-review-2009,larskres-survey-2013} for detailed surveys on tensor decompositions.

Canonical tensor decomposition (also called CP decomposition or CANDECOMP/PARAFAC model) of the eigenvectors of vibrational problems was considered in work by Leclerc and Carrington \cite{lc-rrbpm-2014}.
The authors used rank-reduced block power method (RRBPM).
Each iteration of this method involves matrix-by-vector product, which can be efficiently implemented in tensor formats. The problem is that this method has poor convergence.
Moreover, canonical decomposition is known to have stability issues \cite{desilva-2008}.

The hierarchical RRBPM (H-RRBPM) proposed in~\cite{tc-hrrbpm-2015} by Thomas and Carrington is a significant improvement of the RRBPM.
This method also utilizes sum-of-products representation, but treats strongly coupled coordinates together.
Coupled coordinates are further decomposed hierarchically.

The Multi Configuration Time Dependent Hartree (MCTDH) approach \cite{meyer-book-2009} also uses tensor format, namely the Tucker decomposition. This approach reduces complexity, but suffers from the curse of dimensionality.
This problem was solved in the multilayer MCTDH \cite{wt-mlmctdh-2003} which is similar to the Hierarchical Tucker representation \cite{gras-hsvd-2010}.

We would also like to discuss tensor algorithms for solving  eigenvalue problems developed in mathematical community.
There are two natural directions of solving eigenvalue problems in tensor formats.
One direction is straightforward generalization of iterative methods to tensor arithmetics with rank reduction on each iteration.
For the canonical decomposition power method with shifts was generalized in \cite{beylkin-2002,garcke-mregr-2009} and used in the RRBPM method.
The preconditioned inverse iteration (PINVIT) for tensor formats was considered in  \cite{mach-innereig-2013,ro-hf-2016,ro-crossconv-2015}.
The inverse iteration used in this paper differs from the PINVIT, which is basically preconditioned steepest descent.
Tensor version of the inverse iteration based on iterative solution of arising linear systems was considered in \cite{khst-eigen-2012}.

The PINVIT iteration can explicitly utilize a preconditioner.
The construction of preconditioners in tensor formats for eigenvalue problems was considered in \cite{khst-eigen-2012,khor-prec-2009,tobler-htdmrg-2011,mach-innereig-2013}.
The  approach for a general operator~\cite{mach-innereig-2013}  uses Newton-Schulz iteration in order to find a good preconditioner. However, due to high amount of matrix multiplications, this approach becomes time-consuming.
{In order to construct a preconditioner one can use approximate inverse matrix or approximate solution of linear systems. See \cite{ds-amen-2014,holtz-ALS-DMRG-2012,bnz-pertmr-2014,dc-ttgmres-2013} for solvers in tensor formats.}

The more advanced LOBPCG method was for the first time generalized in \cite{lebedeva-tensornd-2011,tobler-htdmrg-2011}.
We utilize this method and construct preconditioner based on optimization procedure.
The PINVIT method with the proposed preconditioner and the LOBPCG with and without preconditioning were tested in Section~\ref{sec:num}.
Although rank-truncated PINVIT method is fast and is able to accurately calculate small amount of eigenvalues, it fails to converge when a lot of eigenvalues are required.


Alternatively to iterative methods, one can pose an optimization problem -- minimization of the Rayleigh quotient with the constraint on rank.
This approach was recently considered in the matrix product state (MPS) community \cite{pizorn-eigb-2012} and independently  in \cite{dkos-eigb-2014}.
The only disadvantages is that all eigenfunctions are treated in one basis, which leads to large ranks and the method becomes slow (calculation of the acetonitrile took several days).
Nevertheless, this approach becomes very efficient and accurate when small number of eigenpairs are required.

\section{Conclustion and future work} \label{sec:conc}
One of the most interesting missing bits is the theory of the proposed approach. First, why the eigenfunctions can be well-approximated in the TT-format and what are the restrictions on the PES.  Second, the convergence properties of the manifold preconditioner have to be established. These questions will be addressed in future works.

{
From practical point of view, the applicability of the proposed approach for general molecules has to be studied. Currently, it requires the explicit knowledge of sum-of-product representation. Obtaining such a representation is a separate issue, which can be done by using existing methods: POTFIT \cite{jm-potfit-1996} or more general TT-cross approximation approach \cite{ot-ttcross-2010}. If PES has large ranks, coordinate transformation can be helpful \cite{bo-astt-2015}.
}


\section{Acknowledgements} \label{sec:ack}
We would like to thank Prof. Tucker Carrington Jr. and his research group for providing data for the numerical experiment section.
The authors also would like to thank Alexander Novikov for his help in improving the manuscript.
We also thank anonymous referees for their comments and constructive suggestions.

This work has been supported by Russian Science Foundation Grant 14-11-00659.

\appendix
\section{ALS minimization}\label{sec:als}
{
Let us discuss technical details of solving the problem \eqref{eq:alsmin}.
To illustrate the idea let us start from the skeleton decomposition in 2D. In this case $\mX = U V^\top$, where $U,V\in \mathbb{R}^{n\times r}$.
This representation is equivalent to the TT-format in 2D. Indeed, by defining cores $X^{(1)}_\alpha (i_1) = U_{i_1\alpha}$ and $X^{(2)}_\alpha (i_2) = V_{i_2\alpha}$ we get \eqref{eq:tt}.

The ALS procedure starts with fixing one core.
Let us fix $V$ with orthonormal columns (can always be done by QR decomposition) and find the updated $U$ from the minimization problem
\begin{equation}\label{eq:minax}
\min_{U} \| (\mH - \sigma \mathcal{I}) (UV^\top) - \mX_k\|.
\end{equation}
Thanks to well-known property of Kronecker products $(B^T \otimes A)\, \text{vec}(X) = \text{vec}(AXB)$ the problem \eqref{eq:minax} can be equivalently represented as a minimization problem on unknown vector $u\equiv\text{vec}(U)$
$$
\min_{u} \| (H - \sigma I) (V \otimes I_{n})\, u- \text{vec}(\mX_k)\|,
$$
where $\text{vec}(X)$ denotes vectorization of $X$ by reshaping it into a column vector.
Finally we get small linear system
$$
(V^\top \otimes I_{n}) (H - \sigma I)^2 (V \otimes I_{n})\, u = (V^\top \otimes I_{n}) (H - \sigma I)\text{vec}(\mX_k),
$$
where matrix $(V^\top \otimes I_{n}) (H - \sigma I)^2 (V \otimes I_{n})$ is of size $nr \times nr$, while the initial matrix $(H - \sigma I)$ is of size $n^2\times n^2$.
To avoid squared condition number one can formally use the projectional approach
$$
(V^\top \otimes I_{n}) (H - \sigma I) (V \otimes I_{n})\, u = (V^\top \otimes I_{n})\, \text{vec}(\mX_k),
$$
which corresponds to the zero gradient of the energy functional
$$
\mathcal{J}(U, V) = \left<(\mathcal{H} - \sigma \mathcal{I})UV^T,\, UV^T \right> - 2 \left< \mX_k, UV^T\right>.
$$
The problem is that  $(H - \sigma I)$ is not positive definite unless the smallest eigenvalue is required.
Nevertheless, we found that  this approach also works if $UV^T$ is close enough to the eigenvector corresponding to the eigenvalue closest to~$\sigma$. This assumption holds as before running MP II we find a good initial guess for both $\sigma$ and eigenvector, see Section \ref{sec:mp1}.

Note that operator $H$ is also given in the low-rank matrix format (also corresponds to the TT-matrix format in 2D):
$$
H - \sigma I = \sum_{\alpha=1}^R A_\alpha \otimes B_\alpha,
$$
where $A_\alpha, B_\alpha \in \mathbb{R}^{n\times n}$. Hence, the matrix of a local system can be represented as
$$
(V^\top \otimes I_{n}) (H - \sigma I) (V \otimes I_{n}) = \sum_{\alpha=1}^R (V^\top A_\alpha V) \otimes B_\alpha,
$$
and one can do matrix-vector products with~$\mathcal{O}(n (n+r) r R)$ complexity.
Indeed, matrices $(V^\top A_\alpha V)$ are of size $r\times r$ and $B_\alpha$ are $n\times n$. Hence multiplication of $(V^\top A_\alpha V) \otimes B_\alpha$ by a vector  requires $n^2 r + n r^2$ operations.

Similarly to the 2D case in arbitrary dimension we get the following linear system on the vectorized $p$-th core~$x^{(p)} = \text{vec}\,(X^{(p)})$
$$
X_{({\not= p})}^\top (H - \sigma I) X_{(\not= p)} x^{(p)} = X_{({\not= p})}^\top \text{vec}(\mX_k),
$$
where
$$
X_{(\not= p)} = X^{<p} \otimes I_n \otimes X^{>p} \in \mathbb{R}^{n^d\times r_{p-1} n r_p }.
$$
Matrix $X^{<p}$ is a multiplication of first $p-1$ cores reshaped into $n^d \times r_p$ matrix:
$$
	X^{<p} \left(\overline{i_1i_2\dots i_{p-1}},\, :\right) =  X^{(1)}(i_1)\, X^{(2)}(i_2)\, \dots\, X^{(p-1)}(i_{p-1}).
$$
 Matrix $X^{>p}$ is defined by analogy.
Typically additional orthogonalization of $X^{<p}$ and $X^{>p}$ is done \cite{holtz-ALS-DMRG-2012}.

Since operator $(H - \sigma I)$ is given by its TT-representation, matrix-vector multiplication requires $\mathcal{O}(n (dn + r^2) r^2 R^2)$ operations, where $r$ is TT-rank of the tensor $\mX$ and $R$ is the maximum rank of the operator $(H - \sigma I)$. Compared with $\mathcal{O}(n (n+r)r R)$ of the 2D case TT-cores between the first and last dimensions are of size $r\times n \times r$ for TT-tensor and $R\times n \times R$ for TT-matrix, so squared ranks appear.
We note that if dimension and/or mode size are large: $r^2 \ll dn$, then the total complexity is $\mathcal{O}(dn^2 r^2 R^2)$.
}

\newpage
{
\LTcapwidth=\textwidth
\setlength{\tabcolsep}{7pt}
\begin{longtable*}{lcHcHcHcHcHcccc}
\caption{Energy levels and their absolute error (cm$^{-1}$) of acetonitrile for the proposed MP LOBPCG, MP simultaneous inverse iteration and for the H-RRBPM method. Absolute error is presented w.r.t. Smolyak quadrature calculations. H-RRBPM timings are taken from~\cite{tc-hrrbpm-2015}. Gray numbers represent negative error values.}\label{tab:energies} \\
\hline 
 & & \multicolumn{6}{|c|}{H-RRBPM \cite{tc-hrrbpm-2015}} &\multicolumn{2}{c|}{{\footnotesize \textbf{MP LOBPCG}}} & \multicolumn{4}{c|}{\textbf{MP SII}} & Reference  \\
\cline{3-14}
 & &  \multicolumn{2}{|c}{Basis-1} &  \multicolumn{2}{c}{Basis-2} & \multicolumn{2}{c|}{Basis-3} &  \multicolumn{2}{c|}{$r=12$} &  \multicolumn{2}{c}{$r=25$} &  \multicolumn{2}{|c|}{$r=40$}  & Smolyak \cite{ac-smolvibr-2011}\\ 
 & &  \multicolumn{2}{|c}{$6.7$ Mb}    &  \multicolumn{2}{c}{$29$ Mb}  &  \multicolumn{2}{c|}{$139$ Mb} &  \multicolumn{2}{c|}{$9.5$ Mb} &  \multicolumn{2}{c|}{$41$ Mb} &  \multicolumn{2}{c|}{$104$ Mb}  & $1.5$ Gb  \\ 

&  & \multicolumn{2}{|c}{$44$ sec} & \multicolumn{2}{c}{$11$ min} & \multicolumn{2}{c|}{$3.2$ h} &  \multicolumn{2}{c|}{ $17$ min} &  \multicolumn{2}{c}{+$9$ min} &  \multicolumn{2}{|c|}{+$12$ min}  & \\ 
\hline
{} Level & Sym. &       $E$ &    $E - E_\text{ref}$ &       $E$ &    $E - E_\text{ref}$ &       $E$ &     $E - E_\text{ref}$ &            $E$ &     $E - E_\text{ref}$ &        $E$ &     $E - E_\text{ref}$ &            $E - E_\text{ref}$ &  $E$ &      $E_\text{ref}$ \\
\hline
ZPVE &         &  9837.893 &           0.485 &  9837.525 &    0.118 &  9837.429 &           0.022 &  9837.463 &           0.056 &               9837.408 &           0.001 &           \cellcolor[gray]{0.9}{-0.001} &  ${9837.4063}$ &  9837.4073 \\
        $\nu_{11}$ &     $E$ &   361.24 &            0.25 &   361.04 &            0.04 &   361.00 &            0.01 &   361.016 &                        0.03 &                360.991 &                        0.000 &  \cellcolor[gray]{0.9}-0.001 &   360.990 &   360.991 \\
                     &         &   361.24 &            0.25 &   361.06 &            0.07 &   361.00 &            0.01 &   361.029 &                        0.04 &                360.991 &                        0.000 &  \cellcolor[gray]{0.9}-0.001 &   360.990 &   360.991 \\
         $2\nu_{11}$ &     $E$ &   723.60 &            0.42 &   723.41 &            0.23 &   723.22 &            0.04 &   723.274 &                        0.09 &                723.180 &  \cellcolor[gray]{0.9}-0.001 &  \cellcolor[gray]{0.9}-0.001 &   723.180 &   723.181 \\
                     &         &   723.60 &            0.42 &   723.41 &            0.23 &   723.22 &            0.04 &   723.276 &                        0.09 &                723.180 &  \cellcolor[gray]{0.9}-0.001 &  \cellcolor[gray]{0.9}-0.001 &   723.180 &   723.181 \\
         $2\nu_{11}$ &   $A_1$ &   724.27 &            0.44 &   724.05 &            0.23 &   723.87 &            0.04 &   723.919 &                        0.09 &                723.827 &                        0.000 &  \cellcolor[gray]{0.9}-0.001 &   723.826 &   723.827 \\
           $\nu_{4}$ &   $A_1$ &   901.25 &            0.59 &   900.83 &            0.16 &   900.71 &            0.05 &   900.722 &                        0.06 &                900.660 &  \cellcolor[gray]{0.9}-0.002 &  \cellcolor[gray]{0.9}-0.004 &   900.658 &   900.662 \\
           $\nu_{9}$ &     $E$ &  1035.33 &            1.21 &  1034.23 &            0.10 &  1034.19 &            0.07 &  1034.211 &                        0.08 &               1034.127 &                        0.001 &  \cellcolor[gray]{0.9}-0.002 &  1034.124 &  1034.126 \\
                     &         &  1035.34 &            1.21 &  1034.23 &            0.10 &  1034.20 &            0.07 &  1034.241 &                        0.11 &               1034.127 &                        0.001 &  \cellcolor[gray]{0.9}-0.002 &  1034.124 &  1034.126 \\
         $3\nu_{11}$ &   $A_2$ &  1087.23 &            0.67 &  1086.86 &            0.31 &  1086.66 &            0.11 &  1086.720 &                        0.17 &               1086.554 &                        0.000 &  \cellcolor[gray]{0.9}-0.002 &  1086.552 &  1086.554 \\
         $3\nu_{11}$ &   $A_1$ &  1087.23 &            0.67 &  1086.86 &            0.31 &  1086.66 &            0.11 &  1086.734 &                        0.18 &               1086.554 &                        0.000 &  \cellcolor[gray]{0.9}-0.001 &  1086.553 &  1086.554 \\
         $3\nu_{11}$ &     $E$ &  1088.52 &            0.74 &  1088.09 &            0.31 &  1087.88 &            0.11 &  1087.960 &                        0.18 &               1087.776 &                        0.000 &  \cellcolor[gray]{0.9}-0.001 &  1087.775 &  1087.776 \\
                     &         &  1088.52 &            0.75 &  1088.09 &            0.31 &  1087.88 &            0.11 &  1088.005 &                        0.23 &               1087.776 &                        0.000 &  \cellcolor[gray]{0.9}-0.001 &  1087.775 &  1087.776 \\
  $\nu_{4}+\nu_{11}$ &     $E$ &  1260.70 &            0.82 &  1260.04 &            0.16 &  1259.90 &            0.02 &  1259.947 &                        0.07 &               1259.811 &  \cellcolor[gray]{0.9}-0.071 &  \cellcolor[gray]{0.9}-0.073 &  1259.809 &  1259.882 \\
                     &         &  1260.70 &            0.82 &  1260.14 &            0.26 &  1259.90 &            0.02 &  1260.035 &                        0.15 &               1259.811 &  \cellcolor[gray]{0.9}-0.071 &  \cellcolor[gray]{0.9}-0.073 &  1259.809 &  1259.882 \\
           $\nu_{3}$ &   $A_1$ &  1390.32 &            1.34 &  1389.44 &            0.47 &  1389.16 &            0.18 &  1390.003 &                        1.03 &               1389.007 &                        0.034 &  \cellcolor[gray]{0.9}-0.002 &  1388.971 &  1388.973 \\
  $\nu_{9}+\nu_{11}$ &     $E$ &  1396.40 &            1.71 &  1394.90 &            0.21 &  1394.79 &            0.10 &  1395.349 &                        0.66 &               1394.701 &                        0.012 &  \cellcolor[gray]{0.9}-0.007 &  1394.682 &  1394.689 \\
                     &         &  1396.40 &            1.72 &  1394.90 &            0.21 &  1394.80 &            0.11 &  1395.519 &                        0.83 &               1394.709 &                        0.020 &  \cellcolor[gray]{0.9}-0.007 &  1394.682 &  1394.689 \\
  $\nu_{9}+\nu_{11}$ &   $A_2$ &  1396.56 &            1.65 &  1395.10 &            0.19 &  1395.01 &            0.10 &  1395.745 &                        0.84 &               1394.916 &                        0.009 &  \cellcolor[gray]{0.9}-0.007 &  1394.900 &  1394.907 \\
  $\nu_{9}+\nu_{11}$ &   $A_1$ &  1399.39 &            1.70 &  1398.07 &            0.38 &  1397.85 &            0.16 &  1398.570 &                        0.88 &               1397.727 &                        0.040 &  \cellcolor[gray]{0.9}-0.003 &  1397.684 &  1397.687 \\
         $4\nu_{11}$ &     $E$ &  1452.27 &            1.17 &  1451.65 &            0.55 &  1451.35 &            0.25 &  1451.409 &                        0.31 &               1451.095 &  \cellcolor[gray]{0.9}-0.006 &  \cellcolor[gray]{0.9}-0.008 &  1451.093 &  1451.101 \\
                     &         &  1452.27 &            1.17 &  1451.65 &            0.55 &  1451.35 &            0.25 &  1451.501 &                        0.40 &               1451.095 &  \cellcolor[gray]{0.9}-0.006 &  \cellcolor[gray]{0.9}-0.008 &  1451.093 &  1451.101 \\
         $4\nu_{11}$ &     $E$ &  1454.09 &            1.26 &  1453.37 &            0.55 &  1453.08 &            0.25 &  1453.179 &                        0.35 &               1452.821 &  \cellcolor[gray]{0.9}-0.006 &  \cellcolor[gray]{0.9}-0.008 &  1452.819 &  1452.827 \\
                     &         &  1454.24 &            1.41 &  1453.37 &            0.55 &  1453.08 &            0.25 &  1453.231 &                        0.40 &               1452.821 &  \cellcolor[gray]{0.9}-0.006 &  \cellcolor[gray]{0.9}-0.008 &  1452.819 &  1452.827 \\
         $4\nu_{11}$ &   $A_1$ &  1454.85 &            1.45 &  1453.95 &            0.55 &  1453.65 &            0.25 &  1453.739 &                        0.34 &               1453.397 &  \cellcolor[gray]{0.9}-0.006 &  \cellcolor[gray]{0.9}-0.008 &  1453.395 &  1453.403 \\
           $\nu_{7}$ &     $E$ &  1484.53 &            1.30 &  1483.33 &            0.10 &  1483.30 &            0.07 &  1483.455 &                        0.23 &               1483.225 &  \cellcolor[gray]{0.9}-0.004 &  \cellcolor[gray]{0.9}-0.009 &  1483.220 &  1483.229 \\
                     &         &  1484.54 &            1.31 &  1483.34 &            0.11 &  1483.30 &            0.07 &  1483.545 &                        0.32 &               1483.226 &  \cellcolor[gray]{0.9}-0.003 &  \cellcolor[gray]{0.9}-0.008 &  1483.221 &  1483.229 \\
 $\nu_{4}+2\nu_{11}$ &     $E$ &  1621.95 &            1.73 &  1620.86 &            0.64 &  1620.45 &            0.22 &  1620.540 &                        0.32 &               1620.202 &  \cellcolor[gray]{0.9}-0.020 &  \cellcolor[gray]{0.9}-0.024 &  1620.198 &  1620.222 \\
                     &         &  1621.96 &            1.74 &  1620.86 &            0.64 &  1620.45 &            0.22 &  1620.729 &                        0.51 &               1620.202 &  \cellcolor[gray]{0.9}-0.020 &  \cellcolor[gray]{0.9}-0.024 &  1620.198 &  1620.222 \\
 $\nu_{4}+2\nu_{11}$ &   $A_1$ &  1622.56 &            1.79 &  1621.40 &            0.64 &  1620.99 &            0.22 &  1621.281 &                        0.51 &               1620.748 &  \cellcolor[gray]{0.9}-0.019 &  \cellcolor[gray]{0.9}-0.024 &  1620.743 &  1620.767 \\
  $\nu_{3}+\nu_{11}$ &     $E$ &  1751.15 &            1.62 &  1750.18 &            0.65 &  1749.76 &            0.23 &  1750.918 &                        1.39 &               1749.591 &                        0.061 &  \cellcolor[gray]{0.9}-0.005 &  1749.525 &  1749.53\\
                     &         &  1751.15 &            1.62 &  1750.43 &            0.90 &  1749.76 &            0.23 &  1751.229 &                        1.70 &               1749.612 &                        0.082 &  \cellcolor[gray]{0.9}-0.003 &  1749.527 &  1749.53\\
 $\nu_{9}+2\nu_{11}$ &   $A_1$ &  1758.64 &            2.22 &  1756.85 &            0.42 &  1756.60 &            0.18 &  1757.736 &                        1.31 &               1756.475 &                        0.049 &  \cellcolor[gray]{0.9}-0.007 &  1756.419 &  1756.426 \\
 $\nu_{9}+2\nu_{11}$ &   $A_2$ &  1758.65 &            2.22 &  1756.85 &            0.42 &  1756.60 &            0.18 &  1757.756 &                        1.33 &               1756.476 &                        0.050 &  \cellcolor[gray]{0.9}-0.007 &  1756.419 &  1756.426 \\
 $\nu_{9}+2\nu_{11}$ &     $E$ &  1759.26 &            2.13 &  1757.53 &            0.39 &  1757.30 &            0.17 &  1758.023 &                        0.89 &               1757.148 &                        0.015 &  \cellcolor[gray]{0.9}-0.010 &  1757.123 &  1757.133 \\
                     &         &  1759.26 &            2.13 &  1757.53 &            0.40 &  1757.31 &            0.18 &  1758.263 &                        1.13 &               1757.154 &                        0.021 &  \cellcolor[gray]{0.9}-0.009 &  1757.124 &  1757.133 \\
 $\nu_{9}+2\nu_{11}$ &     $E$ &  1762.20 &            2.43 &  1760.37 &            0.59 &  1759.98 &            0.21 &  1761.112 &                        1.34 &               1759.848 &                        0.076 &  \cellcolor[gray]{0.9}-0.004 &  1759.768 &  1759.772 \\
                     &         &  1762.20 &            2.43 &  1760.56 &            0.78 &  1759.98 &            0.21 &  1761.381 &                        1.61 &               1759.859 &                        0.087 &  \cellcolor[gray]{0.9}-0.002 &  1759.770 &  1759.772 \\
          $2\nu_{4}$ &   $A_1$ &  1787.24 &            2.04 &  1785.52 &            0.31 &  1785.37 &            0.16 &  1785.338 &                        0.13 &               1785.128 &  \cellcolor[gray]{0.9}-0.079 &  \cellcolor[gray]{0.9}-0.087 &  1785.120 &  1785.207 \\
         $5\nu_{11}$ &     $E$ &  1818.50 &            1.70 &  1817.82 &            1.02 &  1817.07 &            0.27 &  1817.207 &                        0.41 &               1816.789 &  \cellcolor[gray]{0.9}-0.010 &  \cellcolor[gray]{0.9}-0.012 &  1816.787 &  1816.799 \\
                     &         &  1818.50 &            1.70 &  1817.82 &            1.02 &  1817.07 &            0.27 &  1817.219 &                        0.42 &               1816.789 &  \cellcolor[gray]{0.9}-0.010 &  \cellcolor[gray]{0.9}-0.012 &  1816.787 &  1816.799 \\
         $5\nu_{11}$ &   $A_1$ &  1820.78 &            1.83 &  1819.98 &            1.03 &  1819.22 &            0.27 &  1819.635 &                        0.68 &               1818.943 &  \cellcolor[gray]{0.9}-0.009 &  \cellcolor[gray]{0.9}-0.012 &  1818.940 &  1818.952 \\
         $5\nu_{11}$ &   $A_2$ &  1820.78 &            1.83 &  1819.98 &            1.03 &  1819.22 &            0.27 &  1819.644 &                        0.69 &               1818.944 &  \cellcolor[gray]{0.9}-0.008 &  \cellcolor[gray]{0.9}-0.011 &  1818.941 &  1818.952 \\
         $5\nu_{11}$ &     $E$ &  1821.92 &            1.89 &  1821.07 &            1.04 &  1820.30 &            0.27 &  1820.482 &                        0.45 &               1820.021 &  \cellcolor[gray]{0.9}-0.010 &  \cellcolor[gray]{0.9}-0.014 &  1820.017 &  1820.031 \\
                     &         &  1821.92 &            1.89 &  1821.07 &            1.04 &  1820.30 &            0.27 &  1820.575 &                        0.54 &               1820.021 &  \cellcolor[gray]{0.9}-0.010 &  \cellcolor[gray]{0.9}-0.014 &  1820.017 &  1820.031 \\
  $\nu_{7}+\nu_{11}$ &   $A_2$ &  1845.94 &            1.68 &  1844.45 &            0.19 &  1844.36 &            0.10 &  1845.160 &                        0.90 &               1844.312 &                        0.054 &  \cellcolor[gray]{0.9}-0.008 &  1844.250 &  1844.258 \\
  $\nu_{7}+\nu_{11}$ &     $E$ &  1846.03 &            1.70 &  1844.53 &            0.20 &  1844.43 &            0.10 &  1845.309 &                        0.98 &               1844.387 &                        0.057 &  \cellcolor[gray]{0.9}-0.008 &  1844.322 &  1844.33\\
                     &         &  1846.04 &            1.71 &  1844.53 &            0.20 &  1844.43 &            0.10 &  1845.610 &                        1.28 &               1844.389 &                        0.059 &  \cellcolor[gray]{0.9}-0.008 &  1844.322 &  1844.33\\
  $\nu_{7}+\nu_{11}$ &   $A_1$ &  1846.41 &            1.72 &  1844.89 &            0.20 &  1844.79 &            0.10 &  1845.985 &                        1.30 &               1844.745 &                        0.055 &  \cellcolor[gray]{0.9}-0.009 &  1844.681 &  1844.69\\
   $\nu_{4}+\nu_{9}$ &     $E$ &  1934.56 &            3.01 &  1931.84 &            0.29 &  1931.70 &            0.15 &  1931.849 &                        0.30 &               1931.523 &  \cellcolor[gray]{0.9}-0.024 &  \cellcolor[gray]{0.9}-0.033 &  1931.514 &  1931.547 \\
                     &         &  1934.58 &            3.03 &  1931.84 &            0.29 &  1931.70 &            0.16 &  1931.897 &                        0.35 &               1931.523 &  \cellcolor[gray]{0.9}-0.024 &  \cellcolor[gray]{0.9}-0.032 &  1931.515 &  1931.547 \\
 $\nu_{4}+3\nu_{11}$ &   $A_1$ &  1984.05 &            2.20 &  1983.51 &            1.66 &  1982.14 &            0.29 &  1982.577 &                        0.73 &               1981.821 &  \cellcolor[gray]{0.9}-0.028 &  \cellcolor[gray]{0.9}-0.034 &  1981.815 &  1981.849 \\
 $\nu_{4}+3\nu_{11}$ &   $A_2$ &  1984.05 &            2.20 &  1983.51 &            1.66 &  1982.14 &            0.29 &  1982.783 &                        0.93 &               1981.821 &  \cellcolor[gray]{0.9}-0.029 &  \cellcolor[gray]{0.9}-0.035 &  1981.815 &  1981.85\\
 $\nu_{4}+3\nu_{11}$ &     $E$ &  1985.33 &            2.48 &  1984.51 &            1.65 &  1983.14 &            0.28 &  1983.483 &                        0.63 &               1982.823 &  \cellcolor[gray]{0.9}-0.034 &  \cellcolor[gray]{0.9}-0.041 &  1982.816 &  1982.857 \\
                     &         &  1985.33 &            2.48 &  1984.51 &            1.66 &  1983.14 &            0.29 &  1983.575 &                        0.72 &               1982.823 &  \cellcolor[gray]{0.9}-0.034 &  \cellcolor[gray]{0.9}-0.041 &  1982.816 &  1982.857 \\
          $2\nu_{9}$ &   $A_1$ &  2062.65 &            5.58 &  2058.66 &            1.59 &  2057.17 &            0.10 &  2058.000 &                        0.93 &               2057.071 &                        0.003 &  \cellcolor[gray]{0.9}-0.020 &  2057.048 &  2057.068 \\
          $2\nu_{9}$ &     $E$ &  2069.96 &            4.67 &  2066.78 &            1.49 &  2065.38 &            0.10 &  2065.811 &                        0.53 &               2065.272 &  \cellcolor[gray]{0.9}-0.014 &  \cellcolor[gray]{0.9}-0.019 &  2065.267 &  2065.286 \\
                     &         &  2070.21 &            4.92 &  2066.80 &            1.51 &  2065.38 &            0.10 &  2066.355 &                        1.07 &               2065.297 &                        0.011 &  \cellcolor[gray]{0.9}-0.018 &  2065.268 &  2065.286 \\
 $\nu_{3}+2\nu_{11}$ &     $E$ &  2118.65 &            7.27 &  2112.50 &            1.12 &  2111.75 &            0.37 &  2113.788 &                        2.41 &               2111.511 &                        0.131 &  \cellcolor[gray]{0.9}-0.003 &  2111.377 &  2111.38\\
                     &         &  2118.67 &            7.29 &  2112.50 &            1.12 &  2111.75 &            0.37 &  2113.928 &                        2.55 &               2111.565 &                        0.185 &  \cellcolor[gray]{0.9}-0.001 &  2111.379 &  2111.38\\
 $\nu_{3}+2\nu_{11}$ &   $A_1$ &  2120.01 &            7.71 &  2113.35 &            1.05 &  2112.67 &            0.38 &  2114.145 &                        1.85 &               2112.467 &                        0.170 &  \cellcolor[gray]{0.9}-0.003 &  2112.294 &  2112.297 \\
 $\nu_{9}+3\nu_{11}$ &     $E$ &  2120.94 &            1.61 &  2119.96 &            0.63 &  2119.66 &            0.33 &  2121.397 &                        2.07 &               2119.396 &                        0.069 &  \cellcolor[gray]{0.9}-0.010 &  2119.317 &  2119.327 \\
                     &         &  2121.14 &            1.81 &  2119.96 &            0.63 &  2119.66 &            0.33 &  2122.058 &                        2.73 &               2119.422 &                        0.095 &  \cellcolor[gray]{0.9}-0.010 &  2119.317 &  2119.327 \\
 $\nu_{9}+3\nu_{11}$ &     $E$ &  2121.54 &            1.00 &  2121.10 &            0.56 &  2120.84 &            0.30 &  2122.132 &                        1.59 &               2120.583 &                        0.042 &  \cellcolor[gray]{0.9}-0.012 &  2120.529 &  2120.541 \\
                     &         &  2122.07 &            1.53 &  2121.10 &            0.56 &  2120.85 &            0.31 &  2122.481 &                        1.94 &               2120.586 &                        0.045 &  \cellcolor[gray]{0.9}-0.011 &  2120.530 &  2120.541 \\
 $\nu_{9}+3\nu_{11}$ &   $A_2$ &  2122.07 &            1.16 &  2121.46 &            0.55 &  2121.21 &            0.30 &  2123.442 &                        2.53 &               2120.943 &                        0.033 &  \cellcolor[gray]{0.9}-0.013 &  2120.897 &  2120.91\\
 $\nu_{9}+3\nu_{11}$ &     $E$ &  2123.23 &            0.40 &  2124.17 &            1.34 &  2123.16 &            0.32 &  2124.548 &                        1.71 &               2122.845 &                        0.011 &                        0.001 &  2122.835 &  2122.834 \\
                     &         &  2123.32 &            0.49 &  2124.18 &            1.34 &  2123.16 &            0.33 &  2125.200 &                        2.37 &               2122.997 &                        0.163 &                        0.004 &  2122.838 &  2122.834 \\
 $\nu_{9}+3\nu_{11}$ &   $A_1$ &  2123.59 &            0.29 &  2124.72 &            1.42 &  2123.63 &            0.33 &  2126.877 &                        3.58 &               2123.419 &                        0.118 &  \cellcolor[gray]{0.9}-0.001 &  2123.300 &  2123.301 \\
 $2\nu_{4}+\nu_{11}$ &     $E$ &  2145.31 &            2.70 &  2144.46 &            1.85 &  2142.68 &            0.06 &  2142.585 &  \cellcolor[gray]{0.9}-0.03 &               2142.390 &  \cellcolor[gray]{0.9}-0.224 &  \cellcolor[gray]{0.9}-0.235 &  2142.379 &  2142.614 \\
                     &         &  2145.31 &            2.70 &  2144.60 &            1.99 &  2142.68 &            0.06 &  2143.580 &                        0.97 &               2142.390 &  \cellcolor[gray]{0.9}-0.224 &  \cellcolor[gray]{0.9}-0.235 &  2142.379 &  2142.614 \\
         $6\nu_{11}$ &     $E$ &  2186.73 &            3.09 &  2185.33 &            1.70 &  2184.00 &            0.37 &  2184.164 &                        0.53 &               2183.622 &  \cellcolor[gray]{0.9}-0.013 &  \cellcolor[gray]{0.9}-0.016 &  2183.619 &  2183.635 \\
                     &         &  2186.73 &            3.09 &  2185.33 &            1.70 &  2184.00 &            0.37 &  2184.254 &                        0.62 &               2183.622 &  \cellcolor[gray]{0.9}-0.013 &  \cellcolor[gray]{0.9}-0.016 &  2183.619 &  2183.635 \\
         $6\nu_{11}$ &     $E$ &  2189.99 &            3.86 &  2187.88 &            1.74 &  2186.50 &            0.36 &  2187.113 &                        0.98 &               2186.124 &  \cellcolor[gray]{0.9}-0.014 &  \cellcolor[gray]{0.9}-0.019 &  2186.119 &  2186.138 \\
                     &         &  2190.00 &            3.86 &  2187.88 &            1.74 &  2186.50 &            0.36 &  2187.230 &                        1.09 &               2186.124 &  \cellcolor[gray]{0.9}-0.014 &  \cellcolor[gray]{0.9}-0.019 &  2186.119 &  2186.138 \\
         $6\nu_{11}$ &     $E$ &  2191.88 &            4.24 &  2189.41 &            1.77 &  2188.00 &            0.36 &  2188.710 &                        1.07 &               2187.626 &  \cellcolor[gray]{0.9}-0.016 &  \cellcolor[gray]{0.9}-0.021 &  2187.621 &  2187.642 \\
                     &         &  2192.02 &            4.38 &  2189.42 &            1.78 &  2188.00 &            0.36 &  2188.820 &                        1.18 &               2187.626 &  \cellcolor[gray]{0.9}-0.016 &  \cellcolor[gray]{0.9}-0.021 &  2187.621 &  2187.642 \\
         $6\nu_{11}$ &   $A_1$ &  2192.61 &            4.46 &  2189.93 &            1.78 &  2188.50 &            0.36 &  2189.111 &                        0.97 &               2188.127 &  \cellcolor[gray]{0.9}-0.017 &  \cellcolor[gray]{0.9}-0.022 &  2188.122 &  2188.144 \\
 $\nu_{7}+2\nu_{11}$ &   $A_1$ &  2208.60 &            1.98 &  2207.08 &            0.45 &  2206.79 &            0.17 &  2208.281 &                        1.65 &               2206.740 &                        0.114 &  \cellcolor[gray]{0.9}-0.011 &  2206.615 &  2206.626 \\
 $\nu_{7}+2\nu_{11}$ &   $A_2$ &  2208.61 &            1.98 &  2207.09 &            0.45 &  2206.80 &            0.17 &  2208.938 &                        2.31 &               2206.772 &                        0.139 &  \cellcolor[gray]{0.9}-0.009 &  2206.624 &  2206.633 \\
 $\nu_{7}+2\nu_{11}$ &     $E$ &  2208.73 &            1.97 &  2207.23 &            0.46 &  2206.95 &            0.18 &  2208.948 &                        2.18 &               2206.810 &                        0.044 &                        0.001 &  2206.767 &  2206.766 \\
                     &         &  2208.73 &            1.97 &  2207.23 &            0.46 &  2206.95 &            0.18 &  2209.447 &                        2.68 &               2206.876 &                        0.110 &                        0.002 &  2206.768 &  2206.766 \\
 $\nu_{7}+2\nu_{11}$ &     $E$ &  2209.57 &            2.01 &  2208.01 &            0.45 &  2207.72 &            0.16 &  2209.486 &                        1.93 &               2207.604 &                        0.045 &  \cellcolor[gray]{0.9}-0.010 &  2207.549 &  2207.559 \\
                     &         &  2209.59 &            2.03 &  2208.02 &            0.46 &  2207.72 &            0.16 &  2209.739 &                        2.18 &               2207.615 &                        0.056 &  \cellcolor[gray]{0.9}-0.009 &  2207.550 &  2207.559 \\
\hline 
\end{longtable*}
}




\end{document}